\documentclass[preprint,11pt]{elsarticle}

\usepackage{amssymb}
\usepackage{amsmath}
\usepackage{amsfonts}
\usepackage{amssymb}
\usepackage{indentfirst,latexsym,bm}

\usepackage{tikz}
\usetikzlibrary{calc}

\usepackage{amsthm}
\usepackage{color,xcolor}
\usepackage[all]{xy}
\input{mathrsfs.sty}
\newtheorem{theorem}{Theorem}[section]
\newtheorem{lemma}[theorem]{Lemma}
\newtheorem{corollary}[theorem]{Corollary}
\newtheorem{proposition}[theorem]{Proposition}

\theoremstyle{definition}
\newtheorem{definition}{Definition}[section]

\theoremstyle{definition}

\theoremstyle{remark}
\newtheorem{remark}{Remark}[section]
\theoremstyle{question}

\theoremstyle{problem}

\numberwithin{equation}{section}

\journal{XXX}

\begin{document}

\begin{frontmatter}






\title{The numerical ranges of  the generalized quadratic operators}
\author[shnu]{Kangjian Wu}
\ead{wukjcool@163.com}
\author[shnu]{Qingxiang Xu}
\ead{qingxiang\_xu@126.com}
\address[shnu]{Department of Mathematics, Shanghai Normal University, Shanghai 200234, PR China}
\begin{abstract} We investigate the generalized quadratic operator defined by
$$T =\left(
     \begin{array}{cc}
      a I_H & A \\
c A^* & bI_K\\
       \end{array}
   \right)
,$$  where  $H$ and $K$ are Hilbert spaces, $A:K\to H$ is a bounded linear operator, $I_H$ and $I_K$ denote the identity operators on $H$ and $K$, respectively, and $a,b,c$ are complex numbers. It is shown that $T$ attains its norm if and only if $A$ attains its norm. Furthermore, a complete characterization of the numerical range of $T$ is provided by a new approach.
\end{abstract}

\end{frontmatter}



\section{Introduction}

Throughout this paper, $\mathbb{C}$ is  the complex field,  $\mathbb{N}$ is the set of all positive integers,
 $H$ and $K$ are non-zero complex Hilbert spaces, and $H\oplus K$ represents the Hilbert space equipped with the inner-product defined by
$$\big\langle (x_1,y_1)^T,(x_2,y_2)^T\big\rangle=\langle x_1,x_2\rangle+\langle y_1,y_2\rangle,\quad x_i\in H, y_i\in K, i=1,2.$$
Let $\mathbb{B}(H,K)$ denote the set of all bounded linear operators from $H$ to $K$, abbreviated as  $\mathbb{B}(H)$ when $H=K$. The identity operator on $H$ is denoted by $I_H$, or simply $I$ when no confusion arises. For any $T\in\mathbb{B}(H,K)$, the symbols $T^*$ and $|T|$ represent the adjoint of $T$ and the square root of $T^*T$, respectively. When $T\in\mathbb{B}(H)$, its numerical range $W(T)$ is defined as
$$W(T)=\{\langle Tx,x\rangle: x\in H,\|x\|=1\}.$$
 By the classical Toeplitz-Hausdorff theorem \cite[Theorem 1.7]{WG},  $W(T)$ is convex for every
$ T\in\mathbb{B}(H)$. Furthermore, if $H$ is finite-dimensional, then
$W(T)$ is  closed in $\mathbb{C}$ for any $T\in\mathbb{B}(H)$ \cite[Proposition~1.1]{WG}.

An operator $T\in\mathbb{B}(H\oplus K)$ is called  a generalized quadratic operator \cite[Section~1.1]{LPT} if it admits a block matrix representation  of the form
\begin{equation*}\label{equ:gvTQ}
T =\left(
     \begin{array}{cc}
      a I_H & \lambda A \\
c A^* & bI_K\\
       \end{array}
   \right)\in\mathbb{B}(H\oplus K)
,\end{equation*}
where  $A\in\mathbb{B}(K,H)$, and $a,b,c,\lambda\in\mathbb{C}$.  If $c=0$ or $\lambda=0$, then $T$ reduces to  a quadratic operator in the sense that
$$(T-a I)(T-b I)=0.$$
Therefore, when studying the generalized quadratic operators, we may assume $\lambda\ne 0$. Note that $cA^*$ can be rewritten as
$\frac{c}{\bar{\lambda}}(\lambda A)^*$. Thus, without loss of generality, we may further assume $\lambda=1$. In this case, a generalized quadratic operator takes the form
\begin{align}\label{equ:gvTQ}
T =\left(
     \begin{array}{cc}
      a I_H & A \\
c A^* & bI_K\\
       \end{array}
   \right)\in\mathbb{B}(H\oplus K)
,\end{align}  where  $A\in\mathbb{B}(K,H)$ and $a,b,c\in\mathbb{C}$.

An operator $A \in \mathbb{B}(K,H)$ is said to attain its norm \cite[Section~2.1]{WG} if there exists a unit vector $x \in K $ such that $ \|Ax\|=\|A\|$.
For characterizations of norm attainment, see \cite{BDSS,CN,PP,Ramesh}.
The numerical range of a quadratic operator is completely characterized in \cite[Theorem~2.1]{TW}. In particular, it is
shown there that when $T$ has the form \eqref{equ:gvTQ} with $c=0$, the numerical range $W(T)$ is closed in $\mathbb{C}$ if and only if $T$ attains its norm, which in turn is equivalent to $A$ attaining its norm. Remarkably, reference \cite{TW} has been widely cited in the literature (see e.g.,\,\cite{CTW,Crouzeix,GL,LLPS,LPS,LPT,MT,WPY}). With the exception of the equivalence between the norm attainment of $A$ and $T$, the results of \cite[Theorem~2.1]{TW} have been extended to the case of  the generalized quadratic operators in \cite[Theorem~3.1]{LPT}.  A central aim of this paper is to  complete the theory of \cite[Theorem~3.1]{LPT} by resolving the question of norm attainment. Furthermore, we provide a new approach to deriving the main results in \cite[Theorem~3.1]{LPT}.

The remainder of the paper is organized as follows. Section~\ref{sec:attains norm} addresses the norm attainment of the generalized quadratic operator $T$ defined in \eqref{equ:gvTQ} with its $(1,2)$-entry $A\in\mathbb{B}(K,H)$. Theorem~\ref{thm:morm attainment} establishes that $T$ attains its norm if and only if $A$ attains its norm, thereby generalizing \cite[Lemma~2.2]{TW}. As a consequence of Theorem~\ref{thm:morm attainment}, we establish new equivalent conditions for a quadratic operator to attain its norm, and consequently, a generalization of \cite[Theorem~2.2]{BDSS} is obtained; see Theorem~\ref{thm:equivalent conditions of AN qT} and Corollary~\ref{cor:BDSS} for the details.

Section~\ref{sec:two times two} is devoted to the study of the numerical ranges for a class of $2\times 2$ matrices. Given complex numbers $a, b$, and $c$, a $2\times 2$ matrix $S_d$ is introduced in \eqref{equ:defn of Sd} for each $d\in\mathbb{C}$. It is straightforward to verify that $W(S_d)=W(S_{|d|})$ for all $d\in\mathbb{C}$. Consequently, for every $d>0$, a subset $E_d$ of the complex plane is defined as in \eqref{equ:defn of Ed}. A detailed characterization of $W(S_d)$ and $E_d$ is provided in Lemma~\ref{lem:minor&major lengthes} and Theorem~\ref{thm:summation of E d}, respectively.

Section~\ref{sec:applications} deals with some applications. The numerical range $W(T)$ of a generalized quadratic operator $T$ given by \eqref{equ:gvTQ} is described in Theorem~\ref{thm:nrogqo}. In particular, it is shown that when $|c|=1$,
$a\ne b$, and $c\ne \frac{(a-b)^2}{|a-b|^2}$, the numerical range $W(T)$ is a non-degenerate elliptical disk that is neither open nor closed in $\mathbb{C}$. This result highlights a key difference between the numerical ranges of the generalized quadratic operators and the quadratic operators: according to \cite[Theorem~2.1]{TW}, when $T$ is a quadratic operator, $W(T)$ is always either open or closed in $\mathbb{C}$. Additionally, two propositions concerning the structure of the generalized quadratic operators are provided in this section.

\section{Norm attainments of the  generalized quadratic operators}\label{sec:attains norm}

\begin{lemma}\label{lem:norm of 2*2 operator matrix}{\rm\cite[Lemma~1.6]{FKM}}
Let $T\in \mathbb{B}(H\oplus K)$ be given by
$$
T =\left(
     \begin{array}{cc}
      a I_H & dA \\
c A^* & bI_K\\
       \end{array}
   \right)
,$$  where  $A \in \mathbb{B}(K,H)$ and $a,b,c,d \in \mathbb{C}$. Then\footnote{Note that the term $\|A\|^2$ in the expression for $s$ in \eqref{equ:definition of r and s} was mistakenly written as $\|A\|$ in \cite[Lemma~1.6]{FKM}}
\begin{equation}\label{equ:norm-T}
\|T\| = \frac{1}{2}\left[ (r + s)^{1/2} + (r - s)^{1/2} \right],
\end{equation}
in which
\begin{equation}\label{equ:definition of r and s}
  r= |a|^{2} + |b|^{2} + \|A\|^{2} (|c|^{2}+|d|^{2}), \quad s= 2\Big|ab - cd\|A\|^2\Big|.
\end{equation}

\end{lemma}

\begin{corollary}\label{cor:norm of 2*2 matrix}For every $T=\begin{pmatrix}
         a & d \\
         c & b
       \end{pmatrix}\in M_2(\mathbb{C})$, its norm is given by
\begin{equation}\label{lem:norm of 2*2 matrix}
  \|T\|=\frac{\sqrt{r+s}+\sqrt{r-s}}{2},
\end{equation}
where
\begin{equation*}
  r=|a|^2+|b|^2+|c|^2+|d|^2,\quad s=2|ab-cd|.
\end{equation*}
\end{corollary}
\begin{proof}The conclusion is immediate from Lemma~\ref{lem:norm of 2*2 operator matrix}
by letting $H=K=\mathbb{C}$ and $X=I_K$.
\end{proof}

\begin{lemma}\label{lem:r equiv s}Let $T\in\mathbb{B}(H\oplus K)$ be defined as in \eqref{equ:gvTQ} with its $(1,2)$-entry $A\in\mathbb{B}(K,H)$, and let
   $ \|T\|, r$ and $s$ be given by  \eqref{equ:norm-T} and \eqref{equ:definition of r and s}, respectively. Then
\begin{align}\label{equ:norm-T-squared}&\|T\|^2 = \frac{1}{2}\left(r + \sqrt{r^2 - s^2}\right),\\
\label{equ:rs squared}&r^2 - s^2 = (|a|^2 - |b|^2)^2 + (|c|^2 - 1)^2\|A\|^4+2k\|A\|^2, \\
\label{equ:a plus cb}&|a+\bar{b}c|^2\|A\|^2=\left(\|T\|^2-|b|^2-\|A\|^2\right)\left(\|T\|^2-|a|^2-|c|^2\|A\|^2\right),
\end{align}
where
$$k=\left(|b+\bar{a}c|^2+|a+\bar{b}c|^2\right).$$
\end{lemma}
 \begin{proof}Let $d=1$. By  \eqref{equ:norm-T} and \eqref{equ:definition of r and s}, it follows  the validity of \eqref{equ:norm-T-squared}, and   \begin{align*}&r^2=(|a|^2+|b|^2)^2+(|c|^2+1)^2\|A\|^4+2(|a|^2+|b|^2)(|c|^2+1)\|A\|^2,\\
 &s^2=4\left[|ab|^2+|c|^2\|A\|^4-2\|A\|^2\text{Re}(ab\bar{c})\right].
 \end{align*}
Therefore, the expression for $r^2-s^2$ is given by  \eqref{equ:rs squared}, in which
  \begin{align*}k=&(|a|^2+|b|^2)(|c|^2+1)+4\text{Re}(ab\bar{c})\\
 =&\big[|b|^2+|ac|^2+2\text{Re}(ba\bar{c})\big]+\big[|a|^2+|bc|^2+2\text{Re}(ab\bar{c})\big]\\
 =&|b+\bar{a}c|^2+|a+\bar{b}c|^2.
 \end{align*}

 Using \eqref{equ:norm-T-squared} and the second equation in \eqref{equ:definition of r and s}, we obtain
 \begin{align*}(\|T\|^2 - r)\|T\|^2=&-\frac{1}{4}s^2=-\Big|ab - c\|A\|^2\Big|^2\\
 =&-|ab|^2-|c|^2\|A\|^4+2\text{Re}(ab\bar{c})\|A\|^2.
 \end{align*}
 Let $\gamma$ denote the right-hand side of \eqref{equ:a plus cb}. Then, from the first equation in \eqref{equ:definition of r and s}, we have
    \begin{align*}
    \gamma &= (\|T\|^2 - r)\|T\|^2 + |ab|^2 + (|a|^2 + |bc|^2)\|A\|^2 + |c|^2\|A\|^4 \\
        &= (|a|^2 + |bc|^2)\|A\|^2 + 2\text{Re}(ab\bar{c})\|A\|^2=|a + \bar{b}c|^2\|A\|^2.
    \end{align*}
    This completes the proof.
 \end{proof}

Let $T=\text{diag}(A,B)$ be a diagonal operator, with  $A\in\mathbb{B}(H)$, $B\in\mathbb{B}(K)$, and $\|A\|\le \|B\|$. It is easily seen that $T$ attains its norm if and only if $B$ attains its norm. Moreover, we have the following elementary results concerning norm attainment.
\begin{proposition}\label{prop:attainment-observation}{\rm \cite[Theorem~1.1]{BDSS}} For every $A\in\mathbb{B}(K,H)$, the following statements are equivalent:
\begin{enumerate}
\item[{\rm (i)}]$A$ attains its norm;
\item[{\rm (ii)}]$A^*$ attains its norm;
\item[{\rm (iii)}]$AA^*$ attains its norm;
\item[{\rm (iv)}]$\|A\|^2$ is an eigenvalue of $A^*A$.
\end{enumerate}
\end{proposition}

Now, we provide the main result of this section as follows.
\begin{theorem}\label{thm:morm attainment}Let $T\in\mathbb{B}(H\oplus K)$ be defined as in \eqref{equ:gvTQ} with its $(1,2)$-entry $A\in\mathbb{B}(K,H)$. Then
    the following statements are equivalent:
    \begin{itemize}
        \item [\rm (i)] $T$ attains its norm;
        \item [\rm (ii)] $A$ attains its norm.
    \end{itemize}
\end{theorem}
\begin{proof}If $A=0$, then  both $T$ and $A$ trivially attain their norms. Now, we assume $A\ne 0$.
    A direct computation shows that
        \begin{align*}
        T^*T
        =\left(
           \begin{array}{cc}
            |a|^2 + |c|^2AA^* & (\bar{a} + \bar{c}b)A\\
            (a + c\bar{b})A^* & |b|^2 + A^*A\\
           \end{array}
         \right).
    \end{align*}

    (i) $\Longrightarrow$ (ii).  Suppose $T$ attains its norm. Then $\|T\|^2$ is an eigenvalue of $T^*T$,  so there exists a unit vector $(u,v)^T\in H \oplus K$ such that
    \begin{align*}
        \left(
           \begin{array}{cc}
            |a|^2 + |c|^2AA^* & (\bar{a} + b\bar{c}A\\
            (a + \bar{b}c)A^* & |b|^2 + A^*A\\
           \end{array}
         \right) (u,v)^T= \|T\|^2 (u,v)^T.
    \end{align*}
    This leads to the system:
    \begin{equation}\label{equ:T^*T action on uv}
        \begin{cases}
            (|a|^2 + |c|^2AA^*)u + (\bar{a} + b\bar{c})Av = \|T\|^2 u, \\
            (a + \bar{b}c)A^*u + (|b|^2 + A^*A)v = \|T\|^2 v.
        \end{cases}
    \end{equation}
    Rearranging terms, we obtain:
    \begin{equation}\label{equ:T attaining}
        \begin{cases}
            \left(\|T\|^2 - |a|^2 - |c|^2AA^*\right)u = (\bar{a} + b\bar{c})Av, \\
            \left(\|T\|^2 - |b|^2 - A^*A\right)v = (a + \bar{b}c)A^*u.
        \end{cases}
    \end{equation}
From \eqref{equ:norm-T-squared}, we have $\|T\|^2\ge \frac{r}{2}$, so we consider the following three cases.

    \textbf{Case 1:} $\|T\|^2 = \dfrac{r}{2}$, where $r$ is defined in \eqref{equ:definition of r and s} with $d=1$.
    By \eqref{equ:norm-T-squared}, $r^2-s^2=0$. It follows from \eqref{equ:rs squared} and \eqref{equ:definition of r and s} that
   \begin{align*}&|a|=|b|,\quad |c|=1,\quad a + \bar{b}c=0,\quad r=2(|a|^2+\|A\|^2),
   \end{align*}
   so system \eqref{equ:T attaining} reduces to
   \begin{equation}\label{equ:reduces to trivial case}(\|A\|^2 - AA^*)u = 0, \quad  (\|A\|^2-A^*A)v = 0.\end{equation}
   Since at least one of $u$ or $v$ is non-zero, it follows from Proposition~\ref{prop:attainment-observation} that $A$ attains its norm.

    Suppose that  $\|T\|^2 > \dfrac{r}{2}$. Since $r$ is defined in \eqref{equ:definition of r and s} with $d=1$, we have
    $$\|T\|^2 > |a|^2 + |c|^2\|A\|^2\quad\text{or}\quad \|T\|^2 > |b|^2 + \|A\|^2.$$

    \textbf{Case 2:} $\|T\|^2 > |a|^2 + |c|^2\|A\|^2$. Then the operator $\|T\|^2 - |a|^2 - |c|^2AA^*$ is invertible. From the first equation in \eqref{equ:T attaining}, we get
    \begin{align}\label{equ:v represents u}
        u = (\bar{a} + b\bar{c})\left(\|T\|^2 - |a|^2 - |c|^2AA^*\right)^{-1}Av.
    \end{align}
    This implies $v\ne 0$ (otherwise $u=0$, contradicting that $(u,v)^T$ is a unit vector). Substituting \eqref{equ:v represents u}  into the second equation in \eqref{equ:T attaining} yields
    \begin{align}\label{equ:temp01}
        |a + \bar{b}c|^2(\|T\|^2 - |a|^2 - |c|^2A^*A)^{-1}A^*Av + A^*Av = (\|T\|^2 - |b|^2)v.
    \end{align}
   Since  $$\|T\|^2 - |a|^2 - |c|^2A^*A\ge (\|T\|^2 - |a|^2 - |c|^2\|A\|^2)I,$$ we have
   $$(\|T\|^2 - |a|^2 - |c|^2A^*A)^{-1}\le (\|T\|^2 - |a|^2 - |c|^2\|A\|^2)^{-1} I,$$
   which implies
   $$(\|T\|^2 - |a|^2 - |c|^2A^*A)^{-1}A^*A\le  (\|T\|^2 - |a|^2 - |c|^2\|A\|^2)^{-1}A^*A.$$
   From \eqref{equ:temp01}  and \eqref{equ:a plus cb}, we obtain
   \begin{align*}(\|T\|^2 - |b|^2-\|A\|^2)\|v\|^2=& \langle (\|T\|^2 - |b|^2-\|A\|^2)v,v\rangle\\
   \le&  \langle (\|T\|^2 - |b|^2-A^*A)v,v\rangle\\
   =&|a + \bar{b}c|^2\langle (\|T\|^2 - |a|^2 - |c|^2 A^*A)^{-1}A^*Av,v\rangle\\
   \le &|a + \bar{b}c|^2\langle (\|T\|^2 - |a|^2 - |c|^2\|A\|^2)^{-1}A^*Av,v\rangle\\
   \le &|a + \bar{b}c|^2\langle (\|T\|^2 - |a|^2 - |c|^2\|A\|^2)^{-1}\|A\|^2 v,v\rangle\\
   =&|a + \bar{b}c|^2(\|T\|^2 - |a|^2 - |c|^2\|A\|^2)^{-1}\|A\|^2\|v\|^2\\
   =&(\|T\|^2 - |b|^2-\|A\|^2)\|v\|^2.
   \end{align*}
  Hence, all inequalities become equalities. In particular,
  $$\|A\|^2 \|v\|^2=\|Av\|^2.$$
  Since $v\ne 0$, this implies that $A$ attains its norm.

    \textbf{Case 3:} $\|T\|^2 > |b|^2 + \|A\|^2$. Then
    $\|T\|^2 - |b|^2 - A^*A$ is an invertible operator. From the second equation in \eqref{equ:T attaining}, we obtain
    \begin{align}\label{equ:u represents v}
        v = (a + \bar{b}c)(\|T\|^2 - |b|^2 - A^*A)^{-1}A^*u,
    \end{align}
   which implies $u\ne 0$. Substituting this expression  for $v$ into the first equation in \eqref{equ:T attaining} yields
    \begin{align*}
        |a + \bar{b}c|^2(\|T\|^2 - |b|^2 - AA^*)^{-1}AA^*u = (\|T\|^2 - |a|^2-|c|^2AA^*)u.
    \end{align*}
   Following the same technique as in Case 2, we have  $\|A^*u\|=\|A\|\cdot \|u\|=\|A^*\|\cdot \|u\|$. Therefore,  by Proposition~\ref{prop:attainment-observation},  $A$ attains its norm.

    (ii) $\Longrightarrow$ (i). Assume that $A$ attains its norm. Consider the same cases.

    \textbf{Case 1:} $\|T\|^2 = \frac{r}{2}$.    Choose nonzero vectors $u \in H$ and $v \in K$ such that
    \begin{align*}
        AA^*u = \|A\|^2 u, \quad A^*Av = \|A\|^2 v,\quad \|u\|^2+\|v\|^2=1.
    \end{align*}
  Then  \eqref{equ:reduces to trivial case} is satisfied, so $T$ attains its norm.

    \textbf{Case 2:} $\|T\|^2 > \frac{r}{2}$ and  $\|T\|^2 > |a|^2 + |c|^2\|A\|^2$.
    Choose a nonzero vector $v \in K$ such that $A^*Av = \|A\|^2 v$, and define $u\in H$ via \eqref{equ:v represents u}. Then
    \begin{align*}\left(\|T\|^2 - |a|^2 - |c|^2A^*A\right)v=\left(\|T\|^2 - |a|^2 - |c|^2 \|A\|^2\right)v,
    \end{align*}
    and thus
    $$\left(\|T\|^2 - |a|^2 - |c|^2A^*A\right)^{-1}v=\left(\|T\|^2 - |a|^2 - |c|^2 \|A\|^2\right)^{-1}v.$$
    It follows from \eqref{equ:v represents u}  that
   \begin{align*}A^*u=&(\bar{a} + b\bar{c})\left(\|T\|^2 - |a|^2 - |c|^2A^*A\right)^{-1}A^*Av\\
   =&(\bar{a} + b\bar{c})\left(\|T\|^2 - |a|^2 - |c|^2\|A\|^2\right)^{-1}\|A\|^2v.
   \end{align*}
    By \eqref{equ:a plus cb}, it follows that
    \begin{align*}(a + \bar{b}c)A^*u=&|a + \bar{b}c|^2\left(\|T\|^2 - |a|^2 - |c|^2\|A\|^2\right)^{-1}\|A\|^2v\\
    =&\left(\|T\|^2-|b|^2-\|A\|^2\right)v=\left(\|T\|^2-|b|^2-A^*A\right)v.
    \end{align*}
 This verifies that \eqref{equ:T attaining} holds. Therefore, $(u,v)^T$
 is an eigenvector of $T^*T$ corresponding to $\|T\|^2$. Hence, $T$ attains its norm.

    \textbf{Case 3:} $\|T\|^2 > \frac{r}{2}$ and $\|T\|^2 > |b|^2 + \|A\|^2$. Select a nonzero vector $u \in H$ such that $AA^*u = \|A\|^2 u$, and define $v$ as in \eqref{equ:u represents v}. An argument analogous to that used in Case 2 demonstrates that $T^*T(u,v)^T=\|T\|^2 (u,v)^T$. Therefore, $T$ attains its norm.
\end{proof}

\begin{remark}{\rm A special case of the preceding theorem, which treats only the case $c=0$, appears in \cite[Lemma~2.2]{TW}.
}\end{remark}

As a consequence of Theorem~\ref{thm:morm attainment}, we establish  new equivalent conditions for a quadratic operator to attain its norm.

\begin{theorem}\label{thm:equivalent conditions of AN qT}
Let $Q\in\mathbb{B}(H)$ be a quadratic operator. Then the following statements are equivalent:
\begin{enumerate}
\item[{\rm (i)}] $Q$ attains its norm;
\item[{\rm (ii)}] For any $c,k\in\mathbb{C}$, the operator $Q+cQ^*+kI$ attains its norm;
\item[{\rm (iii)}] There exist $c,k\in\mathbb{C}$ such that the operator $Q+cQ^*+kI$ attains its norm.
\end{enumerate}
\end{theorem}
\begin{proof}By \cite[Theorem~1.1]{TW}, $Q$ is unitarily equivalent to an operator in $\mathbb{B}(H_1\oplus H_2\oplus H_3\oplus H_3)$ of the form
\begin{equation}\label{canonocal form of QT}a_1I_{H_1}\oplus b_1I_{H_2}\oplus \left(
                      \begin{array}{cc}
                        a_1I_{H_3} & A_1 \\
                        0 & b_1I_{H_3} \\
                      \end{array}
                    \right),\end{equation}
where $a_1,b_1\in\mathbb{C}, H_i(1\le i\le 3)$ are Hilbert spaces, and $A_1\in\mathbb{B}(H_3)$ is positive and injective. Since unitarily equivalent operators attain their norms simultaneously,  we may assume that $Q$ is given by \eqref{canonocal form of QT}. For simplicity, we use the same notation $I$ for all identity operators. Thus, for any $c,k\in \mathbb{C}$, we have
$$Q+cQ^*+kI=a_2I\oplus b_2I\oplus T_1,$$
where
\begin{equation}\label{equ:defn of a2 and b2}a_2=a_1+c\overline{a_1}+k,\quad b_2=b_1+c\overline{b_1}+k, \quad T_1=\left(
                                                                                   \begin{array}{cc}
                                                                                     a_2 I & A_1\\
                                                                                     cA_1^* & b_2 I \\
                                                                                   \end{array}
                                                                                 \right).
\end{equation}
Since $\|a_2 I\oplus b_2 I\|\le \|T_1\|$, we conclude that $Q+cQ^*+kI$ attains its norm if and only if $T_1$ attains its norm, which in turn holds if and only if $A_1$ attains its norm (by Theorem~\ref{thm:morm attainment}). Similarly, it follows from \eqref{canonocal form of QT} that $Q$ attains its norm if and only if $A_1$ attains its norm. Therefore, the equivalence is confirmed.
\end{proof}

Since every idempotent  is a quadratic operator, a direct application of Theorem~\ref{thm:equivalent conditions of AN qT} yields the following known result.

\begin{corollary}\label{cor:BDSS}{\rm \cite[Theorem~2.2]{BDSS}} Let $T\in\mathbb{B}(H)$ be an idempotent. Then $T$ attains its norm if and only if the operator $T+T^*-I$ attains its norm.
\end{corollary}

\section{Unions of the numerical ranges of certain $2\times 2$ matrices}\label{sec:two times two}

For an arbitrary subset $E$ of $\mathbb{C}$, we denote by $\overline{E}$, $\text{int}(E)$ and $\partial E$ the closure, interior and boundary of $E$, respectively.
Let $a$, $b$ and $c$ be fixed complex numbers. For any $d\in\mathbb{C}$,  define $S_d\in M_2(\mathbb{C})$ by
\begin{equation}\label{equ:defn of Sd}
S_d =\left(
       \begin{array}{cc}
         a & d \\
         c\bar{d} & b \\
       \end{array}
     \right).
\end{equation}

First, we focus on the characterization  of  the numerical range $W(S_d)$. If $d=0$, then $W(d)$ is clearly a line segment with endpoints $a$ and $b$. For $d\ne 0$, we write $d$ in polar form  $d=|d|e^{i\theta}$, where  $\theta\in [0,2\pi)$. Then, we define  a unitary matrix
$$
U_d=\left(
      \begin{array}{cc}
        1 & 0 \\
        0 & e^{-i\theta}\\
      \end{array}
    \right),
$$
which satisfies
\[
U_d^*S_dU_d = \begin{pmatrix}
    a & |d| \\
    c|d| & b
\end{pmatrix}=S_{|d|}.
\]
It follows that
\begin{equation}\label{equ:radical prop of Sd}W(S_d) =W(S_{|d|}),\quad \forall d\in \mathbb{C}.\end{equation}
Therefore, in the remainder of this section, we assume $d\ge 0$.

\begin{lemma}\label{lem:monotone}For each $d>0$, let $S_d$ be defined by \eqref{equ:defn of Sd}. Then $a,b\in W(S_d)$, and $W(S_{d_1})\subseteq W(S_{d_2})$ whenever $0\le d_1<d_2$.
\end{lemma}
\begin{proof}For any $d>0$, it is clear that
$$ a = \langle S_d(1,0)^T, (1,0)^T\rangle,\quad b = \langle S_d(0,1)^T, (0,1)^T\rangle,$$
so $a,b\in W(S_d)$.

Suppose that $0\le d_1<d_2$. Let $t=\frac{d_1}{d_2}$; then $t\in [0,1)$. Given any $z\in W(S_{d_1})$, there exist $\xi, \eta \in \mathbb{C}$ with $|\xi|^2 + |\eta|^2 = 1$ such that
\begin{align}\label{equ:inner form of z}z=&
        \langle S_{d_1}(\xi,\eta)^T, (\xi,\eta)^T\rangle=a|\xi|^2+d_1(\bar{\xi}\eta+c\xi\bar{\eta})+b|\eta|^2.
\end{align}
Since
$$\langle S_{d_2}(\xi,\eta)^T, (\xi,\eta)^T\rangle=a|\xi|^2+d_2(\bar{\xi}\eta+c\xi\bar{\eta})+b|\eta|^2,$$
and $d_1=td_2$, we have $z=(1-t)z_1+tz_2$, where
$$z_1=a|\xi|^2 + b|\eta|^2,\quad z_2=\langle S_{d_2}(\xi,\eta)^T, (\xi,\eta)^T\rangle.$$
As  $a,b\in W(S_{d_2})$ and $|\xi|^2+|\eta|^2=1$, by the convexity of $W(S_{d_2})$, we conclude that $z_1\in W(S_{d_2})$. This, together with $z_2\in W(S_{d_2})$ and the convexity of $W(S_{d_2})$, yields $z\in W(S_{d_2})$. By the arbitrariness of $z$, we conclude that $W(S_{d_1})\subseteq W(S_{d_2})$.
\end{proof}

The following result is known as the elliptical range theorem. A proof can be found  in \cite[Theorem]{Li} and \cite[Theorem~1.5]{WG}.

\begin{lemma}\label{lem:the numeriacal of matrix01}
Let $A$ be a $2\times2$ matrix with eigenvalues $\lambda_1$ and $\lambda_2$. Then the  numerical range of $A$ is a closed elliptical disk whose foci are $\lambda_1$ and $\lambda_2$, and whose minor axis has length $$\big(\text{tr}(A^*A)-|\lambda_1|^2-|\lambda_2|^2\big)^{\frac12},$$ where $\text{tr}(A^*A)$ denotes the trace of $A^*A$.
\end{lemma}

Given any $d>0$, by Lemma~\ref{lem:the numeriacal of matrix01}, $W(S_d)$ forms an elliptical disk. This includes the degenerate case where the minor axis has length zero; thai is, $W(S_d)$ reduces to a line segment.

\begin{lemma}\label{lem:minor&major lengthes} For each $d>0$, let $S_d$ be defined by \eqref{equ:defn of Sd}. Then $W(S_d)$ is a closed elliptical disk (including the degenerate case) with foci at $\lambda_1(d)$ and $\lambda_2(d)$, and the lengths  of its minor and major axes are given by $2b_d$ and $2a_d$, respectively, where
\begin{align}
    &\lambda_1(d)= \frac{a+b + \sqrt{(a-b)^2+4cd^2}}{2}, \label{equ:lambda1}\\
    &\lambda_2(d)= \frac{a+b - \sqrt{(a-b)^2+4cd^2}}{2}, \label{equ:lambda2}\\
    & b_d = \frac{1}{2\sqrt{2}}\sqrt{
            |a-b|^2 + 2(1 + |c|^2)d^2 - \big|(a-b)^2+4cd^2\big|}, \label{equ:bd} \\
    &a_d = \frac{1}{2\sqrt{2}}\sqrt{
            |a-b|^2 + 2(1 + |c|^2)d^2 +\big|(a-b)^2+4cd^2\big|}.  \label{equ:ad}
\end{align}
\end{lemma}
\begin{proof}Direct computation shows that $\lambda_1(d)$ and $\lambda_2(d)$ are the eigenvalues of $S_d$. By Lemma~\ref{lem:the numeriacal of matrix01} they are the foci of $W(S_d)$. Let $U\in M_2(\mathbb{C})$ be a unitary such that
$$
    U^*S_dU = \left(
                \begin{array}{cc}
                  \lambda_1(d) & \mu_d \\
                  0 & \lambda_2(d) \\
                \end{array}
              \right).$$
Then
$$\mathrm{tr}(S_d^*S_d) = |\lambda_1(d)|^2 + |\lambda_2(d)|^2 + |\mu_d|^2.
$$
On the other hand,
$$
    S_d^*S_d =\left(
                \begin{array}{cc}
                   |a|^2 + |c|^2d^2 & (\bar{a}+b\bar{c})d\\
                   (a+\bar{b}c)d & |b|^2 +d^2\\
                \end{array}
              \right),
     $$
which implies
$$
    \mathrm{tr}(S_d^*S_d) = |a|^2 + |b|^2 + (1 + |c|^2)d^2.
$$
Thus,
    \begin{align*}
        |\mu_d|^2 &= \mathrm{tr}(S_d^*S_d) - (|\lambda_1(d)|^2 + |\lambda_2(d)|^2) \\
        &=|a|^2 + |b|^2 + (1 + |c|^2)d^2- (|\lambda_1(d)|^2 + |\lambda_2(d)|^2).
    \end{align*}
Note that
\begin{align*}\lambda_1(d)+\lambda_2(d)=a+b,\quad \lambda_1(d)-\lambda_2(d)=\sqrt{(a-b)^2+4cd^2}.
\end{align*}
Therefore,
\begin{align*}|\lambda_1(d)|^2 + |\lambda_2(d)|^2=&\frac12\left[\left|\lambda_1(d)+\lambda_2(d)\right|^2+\left|\lambda_1(d)-\lambda_2(d)\right|^2\right]\\
=&\frac12\left[|a+b|^2+\big|(a-b)^2+4cd^2\big|\right].
\end{align*}
Since
$$|a|^2+|b|^2=\frac12\left(|a+b|^2+|a-b|^2\right),
$$
it follows that
\begin{align*}
        |\mu_d|^2 = \frac12 |a-b|^2+ (1 + |c|^2)d^2-\frac12 \left|(a-b)^2+4cd^2\right|.
    \end{align*}
   By Lemma~\ref{lem:the numeriacal of matrix01}, we have  $b_d^2=\frac14 |\mu_d|^2$ and
    \begin{align*}
       a_d^2 &= b_d^2 + \frac{|\lambda_1(d)-\lambda_2(d)|^2}{4}= b_d^2+\frac14\big|(a-b)^2+4cd^2\big|.
    \end{align*}
 The desired conclusion  follows immediately.
 \end{proof}

 \begin{corollary}\label{cor:condition w(s_d) is line segment}
    For each $d>0$, let $S_d$ be defined as in \eqref{equ:defn of Sd}. Then the following statements hold:
 \begin{enumerate}
 \item [{\rm (i)}] If $a=b$, then $W(S_d)$ degenerates to a line segment if and only if $|c|=1$. In this case,
 $W(S_d)$ is the closed line segment with endpoints
 \begin{equation*}\label{newlambda1-012}
    \lambda_1(d)=a+d\sqrt{c},\quad \lambda_2(d)=a-d\sqrt{c}.
 \end{equation*}

 \item [{\rm (ii)}] If $a\ne b$, then  $W(S_d)$ degenerates to a line segment if and only if $$c=\frac{(a-b)^2}{|a-b|^2}.$$ In this case, $W(S_d)$ is the closed line segment with endpoints
    \begin{align}
        \label{newlambda1-03}\lambda_1(d)&=\frac{a+b+(a-b)\sqrt{1+4kd^2}}{2},\\
          \label{newlambda1-04}\lambda_2(d)&=\frac{a+b-(a-b)\sqrt{1+4kd^2}}{2},
    \end{align}
    where $k=\frac{1}{|a-b|^2}$.
 \end{enumerate}
  \end{corollary}
 \begin{proof}We adopt the notation from Lemma~\ref{lem:minor&major lengthes}. Note that
 $W(S_d)$ degenerates to a line segment if and only if $b_d=0$. Observe that
 \begin{align*}\big|(a-b)^2+4cd^2\big|\le |a-b|^2+4|c|d^2\le |a-b|^2+2(1 + |c|^2)d^2.
 \end{align*}
 By \eqref{equ:bd},  we have
 $$b_d=0\Longleftrightarrow \big|(a-b)^2+4cd^2\big|=|a-b|^2+2(1 + |c|^2)d^2.$$
 Therefore, $b_d=0$ if and only if $|c|=1$ and
 \begin{equation}\label{equ:4linear}\big|(a-b)^2+4cd^2\big|=|a-b|^2+4|c|d^2.\end{equation}
If  $a=b$, then the above equality is trivially satisfied.  If  $a\ne b$, we define
$$k=\frac{c}{(a-b)^2}.$$
Then equality \eqref{equ:4linear}  reduces to
$$|1+4kd^2|=1+4|k|d^2,$$
which holds if and only if $k\ge 0$. Since we also have $|c|=1$, this is equivalent to
$$k=\left|\frac{c}{(a-b)^2}\right|=\frac{1}{|a-b|^2}.$$

Now suppose that $W(S_d)$ is a line segment. As shown above,  the endpoints of the line segment are exactly the foci $\lambda_1(d)$ and $\lambda_2(d)$. When $a\ne b$, we have $c=k(a-b)^2$.  Hence, the desired formulas for $\lambda_1(d)$ and $\lambda_2(d)$ follows from
\eqref{equ:lambda1} and \eqref{equ:lambda2}.
 \end{proof}

\begin{remark}{\rm  For each $d>0$, let $S_d$ be defined by \eqref{equ:defn of Sd}. A direct computation shows that
\begin{align*}&S_dS_d^*=\left(
                          \begin{array}{cc}
                            |a|^2+d^2 & (a\bar{c}+\bar{b})d \\
                            (c\bar{a}+b)d & |c|^2 d^2+|b|^2 \\
                          \end{array}
                        \right),\\
&S_d^*S_d=\left(
            \begin{array}{cc}
              |a|^2+|c|^2d^2 & (\bar{a}+\bar{c}b)d \\
              (a+\bar{b}c)d & d^2+|b|^2 \\
            \end{array}
          \right).
\end{align*}
It follows that $S_d$ is normal if and only if $$|c|=1\quad\text{and}\quad c(\bar{a}-\bar{b})=a-b,$$ which in turn  holds if and only if either $|c|=1$ and $a=b$, or $a\ne b$ and $c=\frac{(a-b)^2}{|a-b|^2}$.
Therefore, Corollary~\ref{cor:condition w(s_d) is line segment} implies that $W(S_d)$ degenerates to a line segment if and only if $S_d$ is a normal matrix. This observation is due to \cite[Theorem~3.1]{LPT}.
}\end{remark}

\begin{remark}\label{rem:a&b ne partial point}{\rm
Let $d>0$ and $S_d$ be defined by \eqref{equ:defn of Sd} such that $W(S_d)$ degenerates to a line segment. Direct computation confirms that $a$ and $b$ do not coincide with the endpoints $\lambda_1(d)$ and $\lambda_2(d)$ specified in either case of Corollary~\ref{cor:condition w(s_d) is line segment}.
}\end{remark}

 \begin{corollary}\label{cor:a&bextreme} For each $d>0$,  let $S_d$ be defined by \eqref{equ:defn of Sd} and suppose that $W(S_d)$ is a non-degenerate elliptical disk. Then    $a\in\partial\big(W(S_d)\big)$ if and only if $b\in\partial\big(W(S_d)\big)$, and this occurs if and only if
 $|c|=1$.
 \end{corollary}
 \begin{proof} We adopt the notation from Lemma~\ref{lem:minor&major lengthes}. Since $W(S_d)$ is a non-degenerate elliptical disk, we have
 $$a\in\partial\big(W(S_d)\big)\Longleftrightarrow \big|a-\lambda_1(d)\big|+\big|a-\lambda_2(d)\big|=2a_d.$$
 From the expressions for $\lambda_1(d)$ and $\lambda_2(d)$ given by \eqref{equ:lambda1} and \eqref{equ:lambda2}, we obtain
 \begin{align*}|a-\lambda_1(d)|^2+|a-\lambda_2(d)|^2=&\frac12|a-b|^2+\frac12\big|(a-b)^2+4cd^2\big|,
 \end{align*}
 and
 \begin{align*}\big(a-\lambda_1(d)\big)\big(a-\lambda_2(d)\big)=&a^2-\big(\lambda_1(d)+\lambda_2(d)\big)a+\lambda_1(d)\lambda_2(d)\\
 =&a^2-(a+b)a+ab-cd^2=-cd^2,
 \end{align*}
 which implies
 $$\big|a-\lambda_1(d)\big|\cdot \big|a-\lambda_2(d)\big|=|c|d^2.
 $$
 It follows that
 $$\big(|a-\lambda_1(d)|+|a-\lambda_2(d)|\big)^2 = \tfrac{1}{2}|a-b|^2 + \tfrac{1}{2}\big|(a-b)^2+4cd^2\big| + 2|c|d^2.
 $$
 On the other hand, by  \eqref{equ:ad}  we have
$$4a_d^2=\frac12  |a-b|^2 + (1 + |c|^2)d^2 +\frac12 \big|(a-b)^2+4cd^2\big|.
$$
Thus, when $W(S_d)$ is non-degenerate,
  $$a\in\partial\big(W(S_d)\big)\Longleftrightarrow  2|c|d^2=(1 + |c|^2)d^2\Longleftrightarrow |c|=1.
  $$
 Moreover, since $$\big|b-\lambda_1(d)\big|+\big|b-\lambda_2(d)\big|=\big|a-\lambda_1(d)\big|+\big|a-\lambda_2(d)\big|,$$
 the above analysis shows that when $W(S_d)$ is non-degenerate,  $b\in\partial\big(W(S_d)\big)$ if and only if $|c|=1$.
\end{proof}

Next, for each positive number $d$, we define a subset $E_d$ of $\mathbb{C}$ as follows.

\begin{definition}\label{def:E_d}
    For every $d>0$, let $E_d$ be the subset of $\mathbb{C}$ defined by
\begin{equation}\label{equ:defn of Ed}E_d=\bigcup_{t\in (0,d)} W(S_t),
\end{equation}
where $S_t$ is defined as in \eqref{equ:defn of Sd} for each $t\in (0,d)$.
\end{definition}

To describe $E_d$ in detail, we require the following lemma, which is also of independent interest.
\begin{lemma}\label{lem:convex subset density}
    Let  $E$ be a closed non-degenerate elliptical disk and $F$ be a convex subset of $E$. If $F$ is dense in $E$ and contains the center of $E$, then the interior of $E$ is contained in $F$.
\end{lemma}

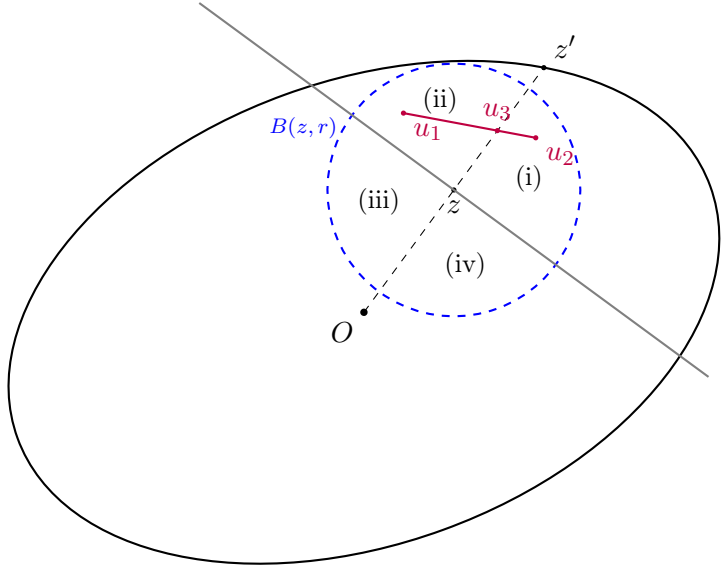
\begin{tikzpicture}[scale=1.4, rotate=20, >=stealth]

    \draw[thick] (0,0) ellipse (3.5 and 2.2);

    \coordinate (O) at (0,0);
    \fill (O) circle (1pt) node[below left] {$O$};

    \coordinate (z) at (1.2,0.8);
    \coordinate (z') at (2.4,1.6);

    \draw[dashed] (O) -- (z');
    \fill (z) circle (0.7pt) node[below] {$z$};
    \fill (z') circle (0.7pt) node[above right] {$z'$};

    \draw[dashed, blue, thick] (z) circle (1.2);
    \node[blue, above right, font=\scriptsize] at (-0.4,1.8) {$B(z,r)$};

    \coordinate (A) at ($(z)!3cm!90:(O)$);
    \coordinate (B) at ($(z)!3cm!-90:(O)$);
    \draw[thick, gray] (A) -- (B);

    \node[font=\small] at ($(z)!0.5!135:(O)$) {(i)};
    \node[font=\small] at ($(z)!0.6!225:(O)$) {(ii)};
    \node[font=\small] at ($(z)!0.5!315:(O)$) {(iii)};
    \node[font=\small] at ($(z)!0.5!45:(O)$) {(iv)};

    \coordinate (u1) at (1.0, 1.65);
    \coordinate (u2) at (2.1, 1);
    \fill[purple] (u1) circle (0.7pt) node[below right] {$u_1$};
    \fill[purple] (u2) circle (0.7pt) node[below right] {$u_2$};
    \draw[purple, thick] (u1) -- (u2);

    \coordinate (u3) at ($(u1)!0.71!(u2)$);
    \fill[purple] (u3) circle (0.7pt) node[above] {$u_3$};

\end{tikzpicture}

\begin{proof}
As shown in the figure, let $O$ be the center of the elliptical disk $E$, and let $z$ be an arbitrary point in $\operatorname{int}(E)$. The ray originating from $O$ and passing through $z$ intersects $\partial E$ at a point $z'$. Since $z$ is an interior point of $E$, there exists $r>0$ such that the open disk $B(z,r)$ centered at $z$ with radius $r$ is entirely contained within $E$.

Through the point $z$, construct a line perpendicular to the ray $Oz$. This line divides the disk $B(z,r)$ into four regions labeled (i), (ii), (iii), and (iv), as depicted. By the density of $F$ in $E$, we can select points  $u_1 \in F$ from region (i) and $u_2 \in F$ from region (ii). The segment joining $u_1$ and $u_2$ intersects the segment $zz'$ at a point $u_3$. Due to the convexity of $F$, the point $u_3$ must belong to $F$.
Since $O \in F$ and $F$ is convex, the entire segment $Ou_3$ is contained in $F$. As $z$ lies between $O$ and $u_3$, it follows that $z \in F$.
By the arbitrariness of $z$, we conclude that $\operatorname{int}(E) \subseteq F$.
\end{proof}

We are now in the position to provide the main result of this section.
\begin{theorem}\label{thm:summation of E d}For every $d>0$, let $S_d$ and $E_d$ be defined by \eqref{equ:defn of Sd} and \eqref{equ:defn of Ed}, respectively. Then $\overline{E_d}=W(S_d)$. Furthermore, the following statements hold:
\begin{enumerate}
\item[\rm (i)] If $|c| \ne 1$,  then $a,b\in \mathrm{int}\big(W(S_d)\big)=E_d$. In this case, $E_d$ is a non-degenerate  open elliptical disk with foci $\lambda_{1}(d)$, $\lambda_{2}(d)$, semi-minor axis $b_d$, and semi-major axis $a_d$ as defined in \eqref{equ:lambda1}--\eqref{equ:ad}.

\item[\rm (ii)] If $a\ne b$ and $c=\frac{(a-b)^2}{|a-b|^2}$, then $E_d$ is the open line segment with endpoints $\lambda_1(d)$ and $\lambda_2(d)$ given by \eqref{newlambda1-03} and \eqref{newlambda1-04}, respectively.

\item[\rm (iii)] If $a = b$ and $|c| = 1$, then $E_d$ is the open line segment with endpoints $a + d\sqrt{c}$ and $a - d\sqrt{c}$.

\item[\rm (iv)] If $a\neq b$, $|c| = 1$, and $c\ne \frac{(a-b)^2}{|a-b|^2}$, then $a,b\in\partial\big(W(S_d)\big)$ and
$$
    E_d = \mathrm{int}\big(W(S_d)\big) \cup \{a, b\},
$$
which is an elliptical disk sharing the same foci and semi-axes as described in part $\rm (i)$. Moreover, in this case $E_d$ is neither open nor closed in $\mathbb{C}$.

\end{enumerate}
\end{theorem}

\begin{proof}
For every $t > 0$, the numerical range $W(S_t)$ is convex.
By Lemma~\ref{lem:monotone}, $W(S_t)$ increases monotonically with $t$.
Hence, Equation \eqref{equ:defn of Ed} implies that $E_d$ is a convex subset of $W(S_d)$.

Let $\{d_n\}_{n=1}^\infty$ be an arbitrary sequence chosen in $(0,d)$ such that $d_n \to d$ as $n\to\infty$. For any
$z\in W(S_d)$, there exists a unit vector $(\xi, \eta)^T$ in $\mathbb{C}^2$ such that
$$z =\big\langle S_d(\xi, \eta)^T, (\xi, \eta)^T \big\rangle= a|\xi|^2 + d(\overline{\xi}\eta + c\xi\overline{\eta}) + b|\eta|^2.$$
Then, for each $n\in\mathbb{N}$, the point  $z_n$, defined by
$$z_n=a|\xi|^2 + d_n(\overline{\xi}\eta + c\xi\overline{\eta}) + b|\eta|^2=\big\langle S_{d_n}(\xi, \eta)^T, (\xi, \eta)^T \big\rangle,$$
satisfies $z_n\in W(S_{d_n})\subseteq E_d$ and $z_n\to z$ as $n\to\infty$.
By the arbitrariness of $z$, we conclude that $W(S_d) \subseteq \overline{E_d}$. On the other hand, since $E_d\subseteq W(S_d)$ and $W(S_d)$ is closed in $\mathbb{C}$, we have $\overline{E_d}\subseteq W(S_d)$. Therefore, $\overline{E_d}=W(S_d)$.

(i). Assume $|c|\ne 1$. By Corollary~\ref{cor:condition w(s_d) is line segment} and Lemma~\ref{lem:minor&major lengthes}, $W(S_d)$ is a non-degenerate elliptical disk with foci $\lambda_{1}(d)$ and $\lambda_{2}(d)$, and with semi-minor axis $b_d$ and semi-major axis $a_d$ as defined in \eqref{equ:lambda1}--\eqref{equ:ad}.

By Lemma~\ref{lem:monotone}, \eqref{equ:defn of Ed}, and Corollary~\ref{cor:a&bextreme}, we have
\begin{equation*}\label{equ:a b in E_d}
    a, b \in \mathrm{int}\big(W(S_d)\big) \cap E_d.
\end{equation*}
Since $E_d$ is convex and the center of $W(S_d)$, given by $(a+b)/2$, lies in $E_d$, and since $E_d$ is dense in $W(S_d)$, it follows from Lemma~\ref{lem:convex subset density} that
\begin{equation*}
    \mathrm{int}\big(W(S_d)\big) \subseteq E_d.
\end{equation*}

It remains to show that $E_d$ contains no boundary points of $W(S_d)$. Suppose, for contradiction, that there exists $z \in \partial \big(W(S_d)\big) \cap E_d$. Then, by \eqref{equ:defn of Ed}, there exist $d_1 \in (0, d)$ and a unit vector $(\xi, \eta)^T \in \mathbb{C}^2$ such that $z$ is given by \eqref{equ:inner form of z}. Consequently, $z$ can be expressed as the convex combination
\begin{align}\label{equ:convex combination for extreme point}
    z = \Big(1 - \frac{d_1}{d}\Big)w_1 + \frac{d_1}{d}w_2,
\end{align}
where
$$w_1=a|\xi|^2 + b|\eta|^2, \quad  w_2=
      a|\xi|^2 + d(\overline{\xi}\eta + c\xi\overline{\eta}) + b|\eta|^2.
$$
As $a,b\in W(S_d)$, $|\xi|^2+|\eta|^2=1$ and
\begin{align*}
    w_2=\big\langle S_d(\xi, \eta)^T, (\xi, \eta)^T \big\rangle,
\end{align*}
we have $w_1,w_2\in  W(S_d)$. Since $W(S_d)$ is a non-degenerate elliptical disk and $z$ is a boundary point, $z$ is an extreme point of $W(S_d)$. Therefore, Equation \eqref{equ:convex combination for extreme point} implies $z = w_1 = w_2$. Hence,
$z = a|\xi|^2 + b|\eta|^2$.
This, together with the fact that $z$ is an extreme point, forces
\begin{equation*}\left\{
                   \begin{array}{ll}
                     z=a, & \hbox{if $|\xi|=1$}, \\
                     z=b, & \hbox{if $|\eta|=1$},\\
                     z=a=b, & \hbox{if $|\xi|<1$ and $|\eta|<1$}.
                   \end{array}
                 \right.
\end{equation*}
However, this contradicts the earlier conclusion that $a, b \in \mathrm{int}\big(W(S_d)\big)$. Therefore, no such $z$ exists, and we conclude that $E_d = \mathrm{int}\big(W(S_d)\big)$.

(ii). Assume that $a\ne b$ and $c=\frac{(a-b)^2}{|a-b|^2}$. For any $t>0$, Corollary~\ref{cor:condition w(s_d) is line segment} tells us that $W(S_t)$ is the closed line segment with endpoints $\lambda_1(t)$ and $\lambda_2(t)$, where $k=\frac{1}{|a-b|^2}$ and
    \begin{align*}
        \lambda_1(t)&=\frac{a+b+(a-b)\sqrt{1+4kt^2}}{2},\\
        \lambda_2(t)&=\frac{a+b-(a-b)\sqrt{1+4kt^2}}{2}.
    \end{align*}
By \eqref{equ:defn of Ed},  $E_d$ is the open line segment with endpoints $\lambda_1(d)$ and $\lambda_2(d)$.

(iii). Assume $a = b$ and $|c| = 1$. By analogy with Part (ii), $E_d$ is the open line segment with endpoints
$a + d\sqrt{c}$ and $a - d\sqrt{c}$.

(iv).  Assume that $a\neq b$, $|c| = 1$, and  $c\ne \frac{(a-b)^2}{|a-b|^2}$. By Lemma~\ref{lem:monotone}, \eqref{equ:defn of Ed}, Corollary~\ref{cor:condition w(s_d) is line segment}, and Corollary~\ref{cor:a&bextreme}, it follows that $W(S_d)$ is a non-degenerate elliptical disk and that
\begin{equation*}\label{equ:a b in E_d}
    a, b \in \partial\big(W(S_d)\big) \cap E_d.
\end{equation*}
By the same reasoning as in Part (i), we obtain
\[
\mathrm{int}\big(W(S_d)\big) \subseteq E_d.
\]
Moreover, a similar argument to that in part (i) shows that no boundary point other than $a$ and $b$ can belong to $E_d$; that is,
\[
E_d \cap \left[\partial \big(W(S_d)\big) \setminus \{a, b\}\right] = \emptyset.
\]
Combining these results, we find
\[
E_d = \mathrm{int}\big(W(S_d)\big) \cup \{a, b\}.
\]
Since $E_d$ contains all interior points of $W(S_d)$ along with exactly two of its boundary points, it follows that $E_d$ is neither an open nor a closed subset of $\mathbb{C}$.
This completes the proof.
\end{proof}

\section{Some applications}\label{sec:applications}

As an application of Theorems~\ref{thm:morm attainment} and \ref{thm:summation of E d}, we first study  the numerical range $W(T)$ of the operator $T\in\mathbb{B}(H\oplus K)$ defined in \eqref{equ:gvTQ}, with its $(1,2)$-entry $A\in\mathbb{B}(K,H)$. Note that if $A=0$, then  $T$ reduces to $aI\oplus bI$,
and thus $W(T)$ is simply the closed line segment joining $a$ and $b$. To exclude this trivial case, we assume in the following theorem that $A\ne 0$.

\begin{theorem}\label{thm:nrogqo}
Let $T \in \mathbb{B}(H \oplus K)$ be defined as in \eqref{equ:gvTQ} with its $(1,2)$-entry $A \in \mathbb{B}(K, H) \setminus \{0\}$. Then
\begin{equation*}\label{equ:either case of W T}
W(T) = W(S_d) \quad \text{or} \quad W(T) = E_d,
\end{equation*}
where $d = \|A\|$, the matrix $S_d$ is given by \eqref{equ:defn of Sd} with its numerical range $W(S_d)$ described in Lemma~\ref{lem:minor&major lengthes}, and $E_d$ is defined in \eqref{equ:defn of Ed} and fully characterized in Theorem~\ref{thm:summation of E d}.
Moreover, the following statements are equivalent:
\begin{enumerate}
\item[\rm (i)] $W(T) = W(S_d)$;
\item[\rm (ii)] $A$ attains its norm;
\item[\rm (iii)] $T$ attains its norm.
\end{enumerate}
\end{theorem}

\begin{proof}
We follow the argument used in the proof of  \cite[Theorem~2.1]{TW}. Let $u$ and $v$ be arbitrary unit vectors in $H$. For any $\xi, \eta \in \mathbb{C}$ satisfying $|\xi|^2 + |\eta|^2 = 1$, we have $\|\xi u\|^2 + \|\eta v\|^2 = 1$ and
\begin{align*}
\left\langle T(\xi u, \eta v)^T, (\xi u, \eta v)^T \right\rangle =
\left\langle
\begin{pmatrix}
a & \langle A v, u \rangle \\
c \langle u, A v \rangle & b
\end{pmatrix}
(\xi, \eta)^T, (\xi, \eta)^T \right\rangle.
\end{align*}
Combining this with \eqref{equ:defn of Sd} and \eqref{equ:radical prop of Sd}, we obtain
\begin{equation}\label{equ:Numerical range of T with big cup}
W(T) = \bigcup_{\substack{u, v \in H \\ \|u\| = \|v\| = 1}} W\big( S_{|\langle A v, u \rangle|} \big).
\end{equation}
Note that
\[
\|A\| = \sup \big\{ |\langle A v, u \rangle| : u, v \in H, \|u\| = \|v\| = 1 \big\},
\]
and as shown in the proof of \cite[Theorem~2.1]{TW}, the above supremum is attained if and only if $A$ attains its norm. Therefore, by Lemma~\ref{lem:monotone} and \eqref{equ:Numerical range of T with big cup}, it follows that either $W(T) = W(S_d)$ or $W(T) = E_d$. Furthermore, $W(T) = W(S_d)$ holds if and only if $A$ attains its norm, which is equivalent to $T$ attaining its norm (see Theorem~\ref{thm:morm attainment}).
This completes the proof.
\end{proof}

Next, we provide two propositions concerning the structure of the generalized quadratic operators.
\begin{proposition}Let $T \in \mathbb{B}(H \oplus K)$ be as defined in \eqref{equ:gvTQ} with its $(1,2)$-entry $A \in \mathbb{B}(K, H)$. If one of the conditions (i)--(iii) in Theorem~\ref{thm:summation of E d} holds, then there exist a quadratic operator $Q\in\mathbb{B}(H\oplus K)$ and a complex number $k$ such that
\begin{equation}\label{equ:equals GQT}T=Q+cQ^*+kI.\end{equation}
\end{proposition}
\begin{proof}Let $Q\in\mathbb{B}(H,K)$ be given by
\begin{equation}\label{constructiong of QT}Q=\left(
                      \begin{array}{cc}
                        a_1I_{H} & A \\
                        0 & b_1I_{K} \\
                      \end{array}
                    \right),\end{equation}
where $a_1,b_1\in\mathbb{C}$. For any $k\in\mathbb{C}$, a direct computation yields
$$Q+cQ^*+kI=\left(
              \begin{array}{cc}
                a_1+c\overline{a_1}+k & A \\
                cA^* & b_1+c\overline{b_1}+k \\
              \end{array}
            \right).$$

\textbf{Case 1:}\ $|c|\ne 0$. Define
$$a_1=\frac{1}{1-|c|^2}(a-c\bar{a}),\quad b_1=\frac{1}{1-|c|^2}(b-c\bar{b}),\quad k=0.$$
Then Equation \eqref{equ:equals GQT} holds.

\textbf{Case 2:}\ $a\ne b$ and $c=\frac{(a-b)^2}{|a-b|^2}$. Set
$$a_1=a-b,\quad b_1=\frac12(a-b),\quad k=2b-a.$$
This satisfies  Equation \eqref{equ:equals GQT}.

\textbf{Case 3:}\  $a = b$ and $|c| = 1$.
Choose
$$a_1=b_1=\frac12 \sqrt{c},\quad k=a-\sqrt{c}.$$
Equation \eqref{equ:equals GQT} follows immediately.
\end{proof}

\begin{proposition}Let $T \in \mathbb{B}(H \oplus K)$ be as defined in \eqref{equ:gvTQ} with its $(1,2)$-entry $A \in \mathbb{B}(K, H) \setminus \{0\}$. If condition (iv)
in Theorem~\ref{thm:summation of E d} holds, then $T$ can not be unitarily equivalent to any operator of the form $Q+cQ^*+kI$, where
$Q$ is a quadratic operator and $k$ is a complex number.
\end{proposition}
\begin{proof} By assumption, we have  $a\neq b$, $|c| = 1$, and $c\ne \frac{(a-b)^2}{|a-b|^2}$. It follows from Theorems~\ref{thm:summation of E d} and \ref{thm:nrogqo} that  $W(T)$ is a non-degenerate elliptical disk.

Suppose, for contradiction, that $T$ is unitarily equivalent  to $Q+cQ^*+kI$ for some quadratic operator $Q$ and complex number $k$.
As in the proof of Theorem~\ref{thm:equivalent conditions of AN qT}, we may assume that $Q$ is given by \eqref{canonocal form of QT}. It follows that
$$Q+cQ^*+kI=a_2I\oplus b_2I\oplus T_1,$$
where $a_2,b_2$ and $T_1$ are derived as in \eqref{equ:defn of a2 and b2}. Since $|c|=1$, we have
$$c(\overline{a_2}-\overline{b_2})=a_2-b_2.$$
Thus, either $a_2=b_2$, or $c=\frac{(a_2-b_2)^2}{|a_2-b_2|^2}$ if $a_2\ne b_2$. Hence, by Theorems~\ref{thm:summation of E d} and \ref{thm:nrogqo}, $W(T_1)$ is a line segment.  Since $W(a_2I\oplus b_2I)\subseteq W(T_1)$, we have
$$W(Q+cQ^*+kI)=W(a_2I\oplus b_2I\oplus T_1)=W(T_1)\ne W(T),$$
which contradicts the assumption that $Q+cQ^*+kI$ and $T$ are unitarily equivalent.
\end{proof}

\end{document}